\documentclass[conference]{IEEEtran}
\usepackage{cite}
\usepackage{amsmath,amssymb,amsfonts}
\usepackage{algorithmic}
\usepackage{graphicx}
\usepackage{textcomp}
\usepackage{xcolor}
\usepackage[version=4]{mhchem}

\def\BibTeX{{\rm B\kern-.05em{\sc i\kern-.025em b}\kern-.08em
    T\kern-.1667em\lower.7ex\hbox{E}\kern-.125emX}}
\begin{document}

\title{A Fully Adaptive Radau Method for the Efficient Solution of Stiff Ordinary Differential Equations at Low Tolerances
}

\author{\IEEEauthorblockN{Shreyas Ekanathan}
\IEEEauthorblockA{Lexington High School \\
Lexington, Massachusetts, USA \\
shreyase39@gmail.com}
\and
\IEEEauthorblockN{Oscar Smith}
\IEEEauthorblockA{JuliaHub \\
Cambridge, Massachusetts, USA \\
oscar.smith@juliahub.com}
\and
\IEEEauthorblockN{Christopher Rackauckas
\IEEEauthorblockA{JuliaHub \\
Massachusetts Institute of Technology \\
Pumas AI inc. \\
Cambridge, Massachusetts, USA \\
chris.rackauckas@juliahub.com}
}
}

\maketitle

\begin{abstract}
Radau IIA methods, specifically the adaptive order Radau method in Fortran due to Hairer, are known to be state-of-the-art for the high-accuracy solution of highly stiff ordinary differential equations (ODEs). However, the traditional implementation was specialized to a specific range of tolerance, in particular only supporting 5th, 9th, and 13th order versions of the tableau and only derived in double precision floating point, thus limiting the ability to be truly general purpose for high fidelity scenarios. To alleviate these constraints, we implement an adaptive-time adaptive-order Radau method which can derive the coefficients for the Radau IIA embedded tableau to any order on the fly to any precision. Additionally, our Julia-based implementation includes many modernizations to improve performance, including improvements to the order adaptation scheme and improved linear algebra integrations. In a head-to-head benchmark against the classic Fortran implementation, we demonstrate our implementation is approximately 2x across a range of stiff ODEs. We benchmark our algorithm against several well-reputed numerical integrators for stiff ODEs and find state-of-the-art performance on several test problems, with a 1.5-times speed-up over common numerical integrators for stiff ODEs when low error tolerance is required. The newly implemented method is distributed in open source software for free usage on stiff ODEs. 
\end{abstract}

\section{Introduction}

High-precision solving of ordinary differential equations (ODEs) is required in many disciplines, from controls and safety scenarios such as identifying potential asteroid collisions \cite{1}, to robustly resolving chemical reactions in astrochemistry which span more time scales than which 64-bit floating point numbers can account for \cite{2}. There are many traditional scenarios where high precision is required as well, for example, in pharmacometrics it is a requirement during inference to solve the model to a much lower tolerance than what's generally used in simulation (tolerances of around $10^{-10}$) in order to improve the robustness of the fitting schemes \cite{3}. This is due to the fact that calculating derivatives of ODE solutions incurs error from the forward pass and induced solver error from the extended equations for the forward/adjoint pass, making higher forward simulation accuracy a requirement in order to achieve accurate gradients and thus a stable optimization scheme \cite{4} \cite{5}.

There is a wealth of high-order schemes for non-stiff ODEs, including the higher order Verner methods \cite{6}, dop853 \cite{7}, high order Gauss-Legendre schemes \cite{8}, and adaptive-order adaptive-time Adams-Bashforth-Moulton schemes like those developed in LSODE \cite{9} and CVODE \cite{10} \cite{11}. Higher order schemes tend to be more efficient at computing solutions to higher accuracy due to the faster convergence with respect to shortened $dt$ \cite{7}. In the SciML open source benchmarks, on non-stiff equations at high accuracy (required errors of roughly 1e-10 and below), schemes above order 5 such as Vern9 typically dominate in terms of performance. 

However, there are many fundamental difficulties for achieving higher-precision solving via high-order methods for stiff equations. Importantly, it is known for many of the standard schemes that an order of 5 is the absolute maximum that can be achieved. For example, with implicit backwards differentiation formulae (BDF) integrators it is known that there is an order barrier where the 6th order method is too unstable for most use cases and all orders above are not zero stable \cite{12}, limiting all standard implementations such as lsoda \cite{13}, CVODE \cite{10} \cite{11}, ode15s \cite{14}, to be an adaptive order scheme where order cannot go above 5 (but is generally bounded below 3 for highly stiff equations for stability region concerns).  Meanwhile, singularly-implicit diagonal Runge-Kutta schemes (SDIRK) and common Rosenbrock-Wanner schemes with single $\gamma$ all have a maximum order of 5 \cite{11}, which is the reason the vast majority of other common stiff ODE solvers are below this order bound such as rodas \cite{15}, the FATODE suite \cite{16}, Rodas5P \cite{17}, and many more. Implicit extrapolation methods could use extra internal stages to overcome this order barrier and achieve arbitrary order, though traditional Fortran implementations such as Seulex and SODEX were never able to achieve the efficiency of the more commonly used ODE solver schemes and the recent Julia advances require using extra compute (exploiting multithreading) in order to achieve state-of-the-art performance \cite{18}.

Thus the only commonly used stiff ODE software for stiff equations which overcomes this order barrier are the Radau IIA schemes, which were popularized by Ernst Hairer's implementation of radau in Fortran \cite{19}. This implementation was a major achievement which uses the fact that Radau IIA schemes of arbitrary order can be derived and are A-B-L stable, giving a platform to build an arbitrary order scheme which is not stability bound for the higher order calculations. However, in order to reduce the computation and due to fundamental programming limitations of the time, the radau scheme pre-derived the schemes for only order 5, 9, and 13 and only for double precision (64-bit) floating point numbers. As such, the scheme tends to perform well for problems defined in double precision which require sufficient accuracy, doing well in the SciML benchmarks for problems which range between 1e-10 to the 64-bit floating point accuracy limit of around 1e-14. In addition, this older limitation does not effectively exploit many of the features of modern hardware, such as pervasive multithreading, single-input multiple data (SIMD) registers, and higher-precision / mixed-precision arithmetic. 

For these reasons, we developed a modern implementation of the adaptive-order adaptive-time Radau IIA schemes in Julia which overcomes the previously mentioned issues and is able to fully exploit modern hardware. The software we have developed is available for use in Julia's OrdinaryDiffEq.jl package, part of the DifferentialEquations.jl interface \cite{20}. Much of the theory detailed in this article is adapted from Hairer's original implementation \cite{7} and \cite{11}. However, our modernization includes a few crucial details, such as the ability to automatically derive order coefficients to arbitrary order to arbitrary user-chosen precision. As such, our adaptive-order scheme can automatically adapt to any order limit chosen by the user, and ultimately generate efficient schemes to any error tolerance provided by the user. Through benchmarks we demonstrate that this scheme achieves state-of-the-art efficiency when highly accurate solutions to stiff ODEs are required.

\section{Background}
\subsection{Fully-Implicit Runge-Kutta Methods}
Consider an initial value problem as follows: 

\begin{equation}
    \frac{\mathrm{d}y}{\mathrm{d}t}=f(t,y) 
\end{equation}
\begin{equation}
    y(t_0)=y_0
\end{equation}
For an $s$ stage Runge-Kutta method with step size $dt$, the approximate solution to $(1)$ at $t$ can be given as:
\begin{equation}
y_{n+1} = y_n + dt \cdot \Sigma_{i=1}^{s}{b_ik_i}
\end{equation}
Where $k_p$ is defined as
\begin{equation}
k_p = f(t_n + c_p \cdot dt, y_n + \Sigma_{i=1}^{p-1}{a_{p,i}k_i})
\end{equation}
for $1 \leq p \leq s$, where the constants $a_{i,j}$ for $1 \leq i, j \leq s$, $b_i$ for $1 \leq i \leq s$, and $c_i$ for $1 \leq i \leq s$ are typically represented by a Butcher tableau \cite{7}, as shown below:
\begin{center}
\renewcommand\arraystretch{1.5}
$\begin{array}
{c|ccccc}
c_1 & a_{11} & a_{12} & \cdots & a_{1,s} \\
c_2 & a_{21} & a_{22} & \cdots & a_{2,s} \\
\vdots & \vdots & \vdots & \vdots & \vdots \\
1 & a_{s,1} & a_{s,2} & \cdots & a_{s,s} \\
\hline 
& b_1 & b_2 & \cdots & b_s 
\end{array}$
\end{center}

Explicit Runge-Kutta methods are methods for which the coefficient matrix $a_{i,j}$ is lower triangular. However, it is a well-known result that such methods have strict stability limits and thus implicit tableaus are required for systems with large condition numbers on the Jacobian, also known as stiffness \cite{11}.

\subsection{Radau IIA Methods}
Radau methods are a class of implicit Runge-Kutta methods that are defined by the Radau quadrature formulas. The tableau is defined by the elements of the Radau polynomial \cite{21}, where $c$ is defined as the roots of the Radau polynomial: 
\begin{equation}
    \frac{\mathrm{d}^{s-1}}{\mathrm{d}x^{s-1}}x^{s-1}(x-1)^s
\end{equation}

We then utilize the so-called "simplifying assumptions" \cite{21} to define two $s$ by $s$ matrices $c_{powers}$ and $c_{q}$ based on $c$, with each element defined as:

\begin{equation}
    c_{powers}[i][j] = c[i]^{j - 1}
\end{equation}

\begin{equation}
    c_q[i][j] = c_{powers}[i][j] \cdot c[i] / j
\end{equation}

Then, the constants $a$ are defined as:
\begin{equation}
    a = c_q \cdot c_{powers}^{-1}
\end{equation} \label{radaueq}
%
%
%
$b$ is simply defined as:
\begin{equation}
    b_i = a_{s,i}
\end{equation}
as required by the order conditions. It can be shown that the $s$ stage Radau method gives an order of $2s-1$ and these methods are A-stable, L-stable, and B-stable for all orders \cite{19}, meaning that all Radau IIA methods generally have higher order than the number of stages and achieve the highest level of stability that is generally provable. It is this high order combined with the with stability that suggests their theoretical potential for highly-accurate integration of stiff systems.

\subsection{Hairer's Fortran Radau Implementation}

The most popular Radau software was developed by Ernst Hairer and published in his landmark paper ``Stiff differential equations solved by Radau methods'' \cite{7}. There, he laid out his implementation of a Radau method that autonomously switches between methods of 5th, 9th, and 13th order \cite{19}. An accompanying implementation in Fortran, simply known as radau, is a standard integrator used in all sorts of platforms from SciPy \cite{22} to OpenModelica \cite{23}. However, there were some limitations to the existing algorithm: 
\begin{enumerate}
    \item The implementation required hard-coded tableau coefficients.
    \item The implementation was limited to methods of 5th, 9th, and 13th order, constraining the flexibility of the method.
\end{enumerate}
A later paper from Martín-Vaquero adds a 17th order method to Hairer's initial implementation, and sees performance boosts at low tolerances \cite{21}. This indicates that adding further adaptive capabilities to the method will continue to improve on performance, motivating our work. 

\section{Implementation of the New Radau Method}

\subsection{Programming Language Choice}

We chose to use Julia to implement our new adaptive Radau method, a contrast to the existing implementation in Fortran. Julia is a high-level programming language which has a simplified interface (garbage collection, built-in array utilities, etc.) like Python or R but has an optimizing compiler and type inference as part of its runtime which allows for its code to be as fast as languages like C or Fortran \cite{24}. This is important for applications such as ODE solving because one of the major bottlenecks to ODE solver performance is the performance of a user's function definition, which when using a high-level programming language like Python will be unoptimized \cite{25}. Additionally, the DifferentialEquations.jl has many advantages over traditional ODE solver libraries, such as support for arbitrary precision arithmetic, automatic differentiation, and automatic sparsity detection which allows for simple user code to achieve maximum performance and precision \cite{20} \cite{26} \cite{27}. Our new adaptive Radau methods make use of all of this infrastructure in tandem, using the arbitrary precision arithmetic to support more accurate convergence with higher orders, automatic differentiation for faster and more accurate Jacobians, and automated sparsity detection as a feature for improving nonlinear solver performance when the underlying $f'$ is sparse. 

\subsection{Arbitrary Order Tableau Generation and Caching}

A notable departure from the previous Radau IIA implementations is the software directly uses the equation \ref{radaueq} for the derivation of the tableaus. The derivation of $a_{i,j}$ is generated using the precision of the provided $u(0)$ initial condition vector and the $c_i$ is generated using the precision of the chosen $(t_0, t_f)$ time span description. The computed tableaus are stored in a dictionary to cache the common tableaus and speed up the setup process for those most frequently used orders, like 5th, 9th, and 13th. If the method chooses to go beyond those orders, it will derive the tableau autonomously, cache the tableau, and continue the numerical integration process. 

\subsection{Embedded Method Generation and Step-Size Adaptivity}

In order to make a Runge-Kutta method adaptive, an alternative reducer $\tilde{b}_i$ vector is typically generated where:
\begin{equation}
\tilde{y}_{n+1} = y_n + dt \cdot \Sigma_{i=1}^{s}{\tilde{b}_ik_i}
\end{equation}
is a lower order approximation. This lower order method is then used to derive an error estimate
\begin{equation}
E = \Vert\tilde{y}_n - y_n\Vert
\end{equation}
which is then used as an approximation for the location truncation error for the adaptive step size controller with a standard predictive controller \cite{7}. Because the embedded Runge-Kutta method to a Radau IIA method is not uniquely defined, we have to make a choice as to the lower order approximation. Following the choice made by Hairer and justified by others \cite{28}, we choose the following linear system of equations that is easily solved for $\tilde{b}$ and gives us a Runge-Kutta method of order $s$:
\begin{equation}
    \Sigma_{j=1}^{s}{\tilde{b}_ic_i^{k-1}} = \frac{1}{k}, k = 1, 2, \cdots, s
\end{equation}

\subsection{Optimizing LU Factorizations via $a$ Transformation}\label{a_transform}

In order to achieve a major speedup over naive implementations, we use a eigen-decomposition trick first exploited in the Hairer radau \cite{7}. To understand the transformation matrix $T$, it is first important to understand the form of $a^{-1}$. Because we only operate on methods with an odd number of stages, $a^{-1}$ always has one real eigenvalue and several complex-conjugate pairs of eigenvalues. For example, when $s=3$, below is $a^{-1}$ and its eigenvectors:
\begin{center}
    $a^{-1} = $         
    $\begin{pmatrix}
         0.197 & - 0.066 & 0.233 \\
         0.394 & 0.292 & -0.042 \\
         0.376 & 0.512 & 0.111 
    \end{pmatrix}$

    $\lambda_1 = 3.6378342527$\\
    $\lambda_2 = 2.681082873627 + 3.0504301992474i$\\
    $\lambda_3 = 2.681082873627 - 3.0504301992474i$\\
\end{center}

The goal now is to find $T$ such that we know the structure of $a^{-1}$ (diagonal, block-diagonal, etc). Naively, we may just take $T$ as the eigenvectors of $a^{-1}$, producing the diagonalization of $a^{-1}$, but this would yield a complex transformation matrix and a complex form of $a^{-1}$, which is suboptimal. Instead, for the below eigenvectors of $a^{-1}$:
\begin{center}
    $\begin{pmatrix}
    0.091 & -0.128+0.027i & -0.128 - 0.027i \\
    0.242 & 0.186 + 0.348i & 0.186-0.348i \\
    0.966 & 0.909 & 0.909
    \end{pmatrix}$
\end{center}
By letting $r$ be the real eigenvector, the complex eigenvectors be represented as $u + vi$ and $u-vi$, where $u$ is a vector holding the real components and $v$ is a vector holding the imaginary components of the complex eigenvectors, we can take a basis ${r, u, v}$, and construct our transformation matrix $T$. For simplicities sake, we also scale $u$ and $v$ to have either a $1$ or a $0$ as their last element, streamlining the computation. Applying this transformation to $a^{-1}$ gives a block diagonal form for $a^{-1}$, with the complex eigenvalues represented as rotation matrices, as shown below:

\begin{center}
    $T = $
    $\begin{pmatrix}
         0.09 & -0.14     &  -0.03 \\
         0.24  &  0.20 &  0.38 \\
         0.96  &  1.0      & 0 \\
    \end{pmatrix}$\\

    $T^{-1}a^{-1} = $
    $\begin{pmatrix}
        3.64 & 0 & 0 \\
        0 & 2.68 & -3.05 \\
        0 & 3.05 & 2.68
    \end{pmatrix}$
\end{center}

Thus, applying our transformation matrix $T$ allows us to have a block-diagonal structure of $a^{-1}$. Thus when moving the Newton solving to the transformed space $z = Tu$, the LU factorizations within the Newton method decouple into $s-1$ $2n \times 2n$ blocks and one $n \times n$ block. Notably, our implementation allows for parallelizing these LU factorizations which can be an optimization when the blocks are too small for BLAS-based multithreading to be effective. However, this transformation is a major performance optimization even without parallelization as the decoupling of the $\mathcal{O}((sn)^3)$ cost to LU factorize the $sn \times sn$ Newton matrix is transformed into $\mathcal{O}(sn^3)$ as each block is done independently. This optimization is thus crucial to our implementation which can achieve high order and thus large stage sizes, effectively decreasing the cost by orders of magnitude for higher fidelity integrations.

\subsection{Order Adaptivity}
In addition to step size adaptivity, our method also builds on the order adaptivity in existing methods by implementing a novel order selection algorithm. Where the Hairer implementation uses the value of the divergence to determine whether to raise or lower the order \cite{19}, our algorithm chooses a different approach. However, the general scheme is similar. The goal of the order adaptivity is to choose the highest order to get the largest time steps for which the Newton method can converge, and then the time step is chosen to be as large as possible in order to keep the error below tolerance. The time-steps of a Radau method uses several Newton iterations while converging on the solution to that timestep, where more iterations implies that the Newton solver is becoming unstable (and thus smaller time steps are required to get a better initial guess) and each iteration requires a fixed cost. Notably, this cost is only linear with respect to the order because of the eigenvector transformation of Section \ref{a_transform}, and thus the cost is both linear with respect to the chosen order and with respect to the number of Newton iterations. Thus we use a metric based on the number of iterations that it takes the method to converge. If the method is taking several iterations to converge, this implies that the Newton solver is at the edge of stability and that a lower order method can be more efficient. 

Technically, we have a heuristic \verb|hist_iter| that tracks the historical number of iterations for out method. Then, we combine the current and historical number of iterations by calculating a parameter $\kappa = \verb|hist_iter| * 0.8 + \verb|iter| * 0.2$. This weighted average allows any given time-step to not dramatically skew the performance. After calculating this parameter, we then perform the following comparisons:
\begin{enumerate}
    \item If $\kappa < 2.75$, then the method is converging sufficiently fast. Raise the order of the method by 4.
    \item Else if $\kappa > 8$, then the method is converging too slowly. Decrease the order of the method by 4.
    \item Otherwise, continue at the same order. 
\end{enumerate}
Note that increments of 4 are used so that every chosen order has the eigenvector structure required in Section \ref{a_transform}. Doing this approach allows us to combine step-size and order adaptivity. When we are taking small timesteps and need less iterations to converge, the algorithm chooses a higher order method, and vice versa. 

\section{Results}
To analyze the results of the method, we use work-precision diagrams, which show the time it takes a given solver to reach a certain degree of accuracy. Lower times and higher accuracies on the diagram indicate better performance. Accuracy is compared against solves of a given problem at tolerances of $10^{-14}$.

\subsection{Competing algorithms}
The algorithms against which we benchmark our adaptive Radau method are described below:
\begin{enumerate}
    \item Hairer's Radau: The existing implementation of Radau methods in Fortran.
    \item BDF Methods: Another set of very common numerical integrators for stiff ODEs. They are multi-step methods that use interpolating polynomials to derive future time-steps. They initially struggled with stability issues, but later improvements have created very powerful methods that perform well on stiff ODEs \cite{29}. 
    \item SUNDIALS BDF Methods: A common C++ implementation of BDF methods. Chosen as a common ``standard candle'' of stiff ODE solvers, as CVODE is a standard method which is also very similar to other methods such as VODE and LSODE.
    \item Rosenbrock methods: A set of numerical integrators that perform only one iteration of Newton's method \cite{29}. Solving the linear system can be expensive in these methods, but when applied correctly, they are among the best numerical integrators for stiff ODEs \cite{12}\cite{29}. 
\end{enumerate}

BDF and Rosebrock methods cannot achieve the same order of accuracy that Radau methods can. Therefore, while BDF and Rosenbrock methods tend to dominate solving stiff ODEs at high tolerances, Radau methods typically dominate the solution to low tolerances. 

\subsection{Work-Precision Diagrams}
We test our adaptive Radau method against several methods belonging to the classes listed above and plot work-precision diagrams. The test problems,  from Julia's SciML benchmarks, are listed below:

\begin{itemize}
    \item Oregonator, 3 ODEs \cite{30} \cite{31}.
    \item Robertson, 3 ODEs \cite{32}.
    \item Hires, 8 ODEs \cite{33}.
    \item Pollution, 20 ODEs \cite{34}. 
\end{itemize}

In each of the diagrams, the algorithm we implemented is labeled "AdaptiveRadau". We run the tests with the possible orders for the method as 5, 9, 13, 17, 21, or 25. 

\begin{figure}[bp!]
\includegraphics[width=\linewidth]{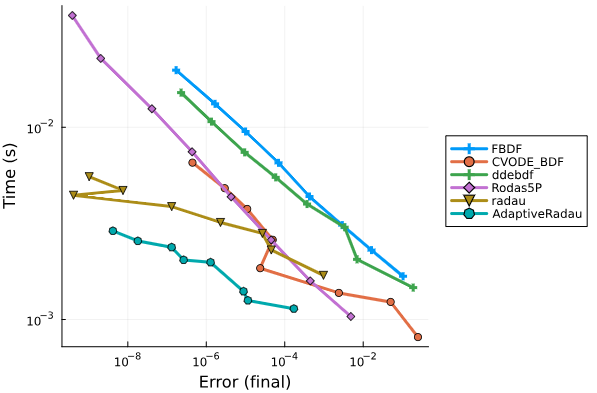}
\caption{Comparison between different Runge-Kutta methods on the Oregonator problem}
\label{fig}
\end{figure}

\begin{figure}[htbp]
\includegraphics[width=\linewidth]{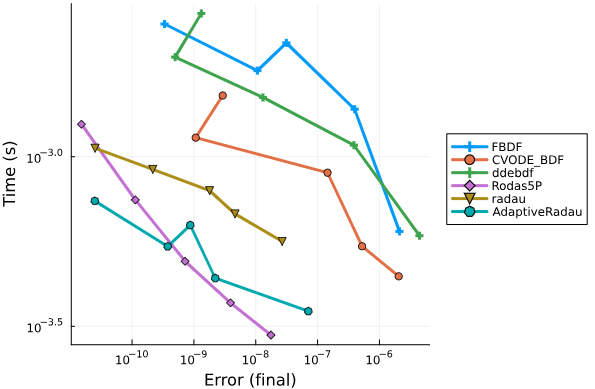}
\caption{Comparison between different Runge-Kutta methods on the Robertson problem}
\label{fig}
\end{figure}

\begin{figure}[htbp]
\includegraphics[width=\linewidth]{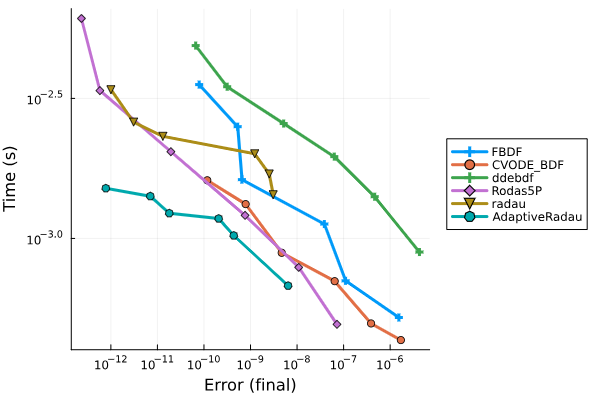}
\caption{Comparison between different Runge-Kutta methods on the Hires problem}
\label{fig}
\end{figure}

\begin{figure}[htbp]
\includegraphics[width=\linewidth]{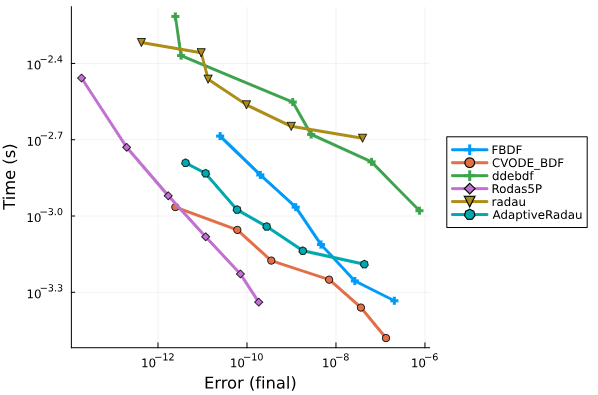}
\caption{Comparison between different Runge-Kutta methods on the Pollution problem}
\label{fig}
\end{figure}

The comparison between only our AdaptiveRadau algorithm and the existing implementation of radau in Fortran is shown in the appendix, where results are consistently better by an average of 2-times across all test problems. 

Furthermore, across all algorithms as a whole, we find that the AdaptiveRadau algorithm we implement achieves state-of-the-art performance on Hires, Oregonator, and Robertson problems. Even when the algorithm does not perform as well, like on the Pollution problem, we still far outperform the existing implementation of Radau. 

We also consider benchmarks of our algorithm beyond the Float64 domain. In the following, we compare our AdaptiveRadau algorithm against one of the best stiff numerical integrators in Julia, Rodas5P, a newer high-order Rosenbrock method \cite{17}. Both solutions are compared with a reference solution computed to tolerances of $10^{-24}$. The problem and solutions are computed with BigFloat precision. Notably, the difference in performance dramatically increases as the required tolerance increases, showcases the special high performance of this new integrator in this regime.

\begin{figure}[htbp]
\includegraphics[width=\linewidth]{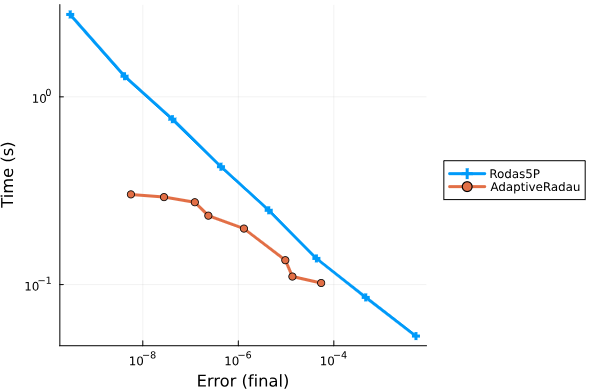}
\caption{Comparison between AdaptiveRadau and Rodas5P on the Oregonator problem at BigFloat precision}
\label{fig}
\end{figure}

\begin{figure}[htbp]
\includegraphics[width=\linewidth]{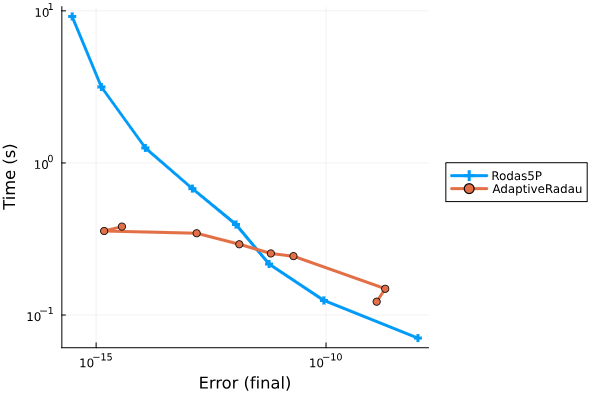}
\caption{Comparison between AdaptiveRadau and Rodas5P on the Hires problem at BigFloat precision}
\label{fig}
\end{figure}

\section{Discussion}

Radau methods have dominated the solution of stiff ODEs at low tolerances since their development because they can achieve much higher orders~\cite{29}. In our work, we see that our algorithm continues to dominate that area. Our methods often are the fastest method to converge on the solution at low error, as seen in the work-precision diagrams above. 

Radau methods had not been changed in decades, since Hairer's implementation appeared to be the gold standard for high-precision stiff ODE solving. However, we now see that the standard has been raised by our algorithm. At high tolerances, we don't outperform the best algorithms for solving stiff ODEs, but neither does Radau. Our implementation does compete with the Rosenbrock methods at these higher tolerances but is often outperformed. However, at low tolerances, the diagrams show a different story. Our algorithm is among the best, if not the best, at solving these stiff ODEs with high precision. 

Furthermore, we show that the mixed-precision capabilities of our algorithm enhance the state-of-the-art. The comparison between Rodas and Radau on the BigFloat problem setups highlights the value of a method that can change orders. The relatively flat curve shown in our algorithm is starkly contrasted with the steeply sloping Rodas5P curve. By adjusting the order to the time steps taken, we optimize our solve. The diagrams show that we are extremely efficient at low tolerances, building on the success we had in the prior diagrams. 

\section*{Acknowledgment}

S.E. thanks the MIT PRIMES Computer Science program for connecting him with his mentor and giving him the chance to work on this research, and Dr. Chris Rackauckas and Mr. Oscar Smith for their continued support throughout the process.

\newpage

\section{Appendix}

\subsection{Direct Comparisons}
Below are comparisons between only the AdaptiveRadau algorithm we implement and the existing implementation of Radau by Hairer in Fortran. The diagrams show that our algorithm is consistently faster than the existing implementation, with speed-ups ranging from roughly 1.5x to 2x.

\begin{figure}[htbp]
\includegraphics[width=\linewidth]{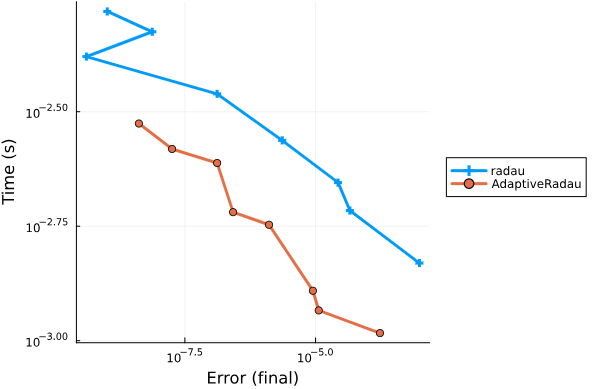}
\caption{Comparison between our AdaptiveRadau method and the existing Radau method on the Oregonator problem}
\label{fig}
\end{figure}

\begin{figure}[htbp]
\includegraphics[width=\linewidth]{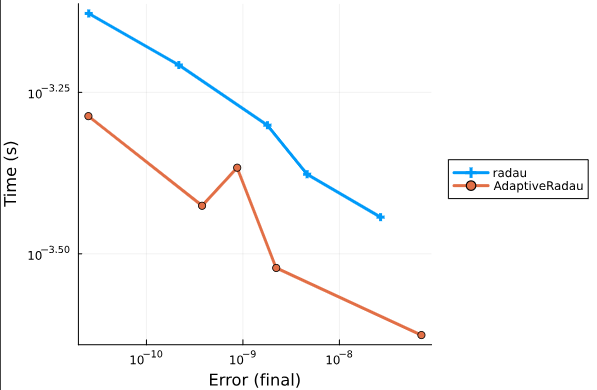}
\caption{Comparison between our AdaptiveRadau method and the existing Radau method on the Robertson problem}
\label{fig}
\end{figure}

\begin{figure}[htbp]
\includegraphics[width=\linewidth]{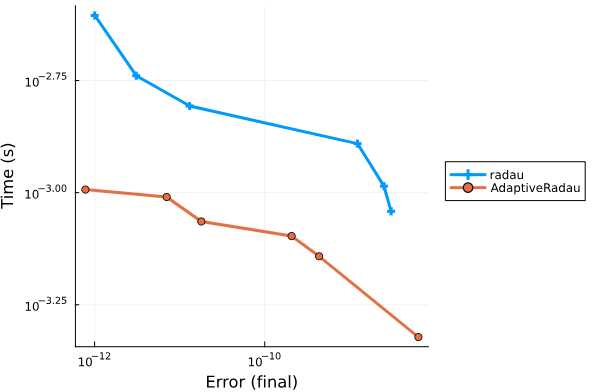}
\caption{Comparison between our AdaptiveRadau method and the existing Radau method on the Hires problem}
\label{fig}
\end{figure}

\begin{figure}[htbp]
\includegraphics[width=\linewidth]{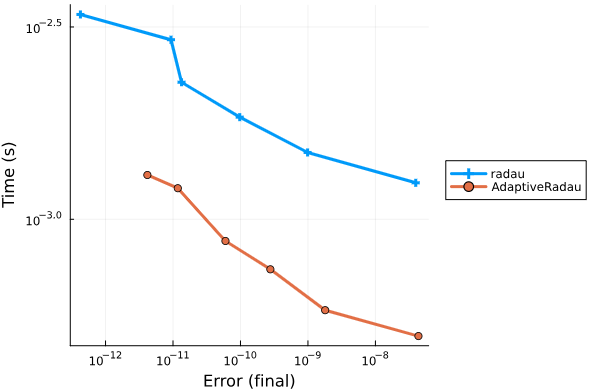}
\caption{Comparison between our AdaptiveRadau method and the existing Radau method on the Pollution problem}
\label{fig}
\end{figure}
\newpage
\subsection{Models and Benchmarks}
Below are the setups for each test problem:

The Oregonator problem is defined as:

\begin{equation}
    \frac{dy_1}{dt} = -k_1(y_2+y_1(1-k_2y_1-y_2))
\end{equation}

\begin{equation}
    \frac{dy_2}{dt}=\frac{y_3-(1+y_1)y_2}{k_1}
\end{equation}

\begin{equation}
    \frac{dy_3}{dt} = k_3(y_1-y_3)
\end{equation}

The initial conditions are $y = [1.0, 2.0, 3.0]$ and $k = (77.27, 8.375 \cdot 10^{-3}, 0.161)$. The time span for integration is $t = (0.0 s, 30 s)$ \cite{30}\cite{31}. 

The comparison with several algorithms is done with relative tolerances in the range $(10^{-5}, 10^{-12})$ and corresponding absolute tolerances in the range $(10^{-7}, 10^{-14})$ against a reference solution with relative and absolute tolerances $10^{-14}$. 

The comparison between only AdaptiveRadau and Rodas5P is done in BigFloat precision with relative tolerances in the range $(10^{-5}, 10^{-12})$ and corresponding absolute tolerances in the range $(10^{-9}, 10^{-16})$ against a reference solution with relative and absolute tolerances of $10^{-24}$.

The Robertson problem is defined as:

\begin{equation}
    \frac{dy_1}{dt} = -k_1y_1+k_3y_2y_3
\end{equation}

\begin{equation}
    \frac{dy_2}{dt}=k_1y_1-k_2y_2^2-k_2y_2y_3
\end{equation}

\begin{equation}
    \frac{dy_3}{dt} = k_2y_2^2
\end{equation}

The initial conditions are $y = [1.0, 0.0, 0.0]$ and $k = (0.04, 3*10^7, 1*10^4)$. The time span for integration is $t = (0.0 s, 1*10^5 s)$ \cite{32}. 

The comparison with several algorithms is done with relative tolerances in the range $(10^{-4}, 10^{-8})$ and corresponding absolute tolerances in the range $(10^{-9}, 10^{-13})$ against a reference solution with relative and absolute tolerances $10^{-14}$. 

The comparison between only AdaptiveRadau and Rodas5P is done in BigFloat precision with relative tolerances in the range $(10^{-5}, 10^{-12})$ and corresponding absolute tolerances in the range $(10^{-9}, 10^{-16})$ against a reference solution with relative and absolute tolerances of $10^{-24}$.

The Hires problem is defined as:

\begin{equation}
    \frac{dy_1}{dt} = -1.71y_1 + 0.43y_2 + 8.32y_3+0.0007
\end{equation}

\begin{equation}
    \frac{dy_2}{dt}=1.71y_1 - 8.75y_2
\end{equation}

\begin{equation}
    \frac{dy_3}{dt} = -10.03y_3 + 0.43y_4 + 0.35y_5
\end{equation}

\begin{equation}
    \frac{dy_4}{dt} =8.32y_2 + 1.71y_3 - 1.12y_4
\end{equation}

\begin{equation}
    \frac{dy_5}{dt}= -1.745y_5 + 0.43y_6 + 0.43y_7
\end{equation}

\begin{equation}
    \frac{dy_6}{dt} = -280.0y_6y_8 + 0.69y_4 + 1.71y_5 - 0.43y_6 + 0.69y_7
\end{equation}

\begin{equation}
    \frac{dy_7}{dt} = 280.0y_6y_8 - 1.81y_7
\end{equation}

\begin{equation}
    \frac{dy_8}{dt} = -280.0y_6y_8 + 1.81y_7
\end{equation}

The initial conditions are $y = [1.0, 0.0, 0.0, 0.0, 0.0, 0.0, 0.0, 0.0057]$. The time span for integration is $t = (0.0 s, 321.8122 s)$ \cite{33}. 

The comparison with several algorithms is done with relative tolerances in the range $(10^{-5}, 10^{-10})$ and corresponding absolute tolerances in the range $(10^{-7}, 10^{-12})$ against a reference solution with relative and absolute tolerances $10^{-14}$. 

The Pollution problem is defined as:

\begin{equation}
    \frac{dy_1}{dt} = -k_1y_1-k_{10}y_{11}k_1-k_{14}y_1y_6-k_{23}y_1y_4
\end{equation}

\begin{equation}
    -k_{24}y_{19}y_1+k_2y_2y_4+k_3y_5y_2+k_9y_{11}y_2
\end{equation}

\begin{equation}
    +k_{11}y_{13}+k_{12}y_{10}y_{2}+k_{22}y_{19}+k_{25}y_{20}
\end{equation}

\begin{equation}
    \frac{dy_2}{dt}=-k_2y_2y_4-k_3y_5y_2-k_9y_{11}y_2-k_{12}y_{10}y_2
\end{equation}

\begin{equation}
    +k_1y_1+k_{21}y_{19}    
\end{equation}

\begin{equation}
    \frac{dy_3}{dt} = -k_{15}y_3+k_1y_1+k_{17}y_4+k_{19}y_{16}+k_{22}y_{19}
\end{equation}

\begin{equation}
    \frac{dy_4}{dt} =-k_2y_2y_4-k_{16}y_4-k_{17}y_4-k_{23}y_1y_4+k_{15}y_3
\end{equation}

\begin{equation}
    \frac{dy_5}{dt}= -k_3y_5y_2+2k_4y_7+k_6y_7y_6+k_7y_9
\end{equation}

\begin{equation}
    +k_{13}y_{14}+k_{20}y_{17}y_6
\end{equation}

\begin{equation}
    \frac{dy_6}{dt} = -k_6y_7y_6-k_8y_9y_6-k_{14}y_1y_6
\end{equation}

\begin{equation}
    -k_{20}y_{17}y_6+k_3y_5y_2+2k_{18}y_{16}
\end{equation}

\begin{equation}
    \frac{dy_7}{dt} = -k_4y_7-k_5y_7-k_6y_7y_6+k_{13}y_{14}
\end{equation}

\begin{equation}
    \frac{dy_8}{dt} = k_4y_7+k_5y_7-k+k_6y_7y_6+k_7y_9
\end{equation}

\begin{equation}
    \frac{dy_9}{dt} = -k_7y_9-k_8y_9y_6
\end{equation}

\begin{equation}
    \frac{dy_{10}}{dt} = -k_{12}y_{10}y_2+k_7y_9+k_9y_{11}y_2
\end{equation}

\begin{equation}
    \frac{dy_{11}}{dt} = -k_9y_{11}y_2-k_{10}y_{11}+k_8y_9y_6+k_{11}y_{13}
\end{equation}

\begin{equation}
    \frac{dy_{12}}{dt} = k_9y_{11}y_2
\end{equation}

\begin{equation}
    \frac{dy_{13}}{dt} = -k_{11}y_{13}+k_{10}y_{11}y_1
\end{equation}

\begin{equation}
    \frac{dy_{14}}{dt} = -k_{13}y_{14}+k_{12}y_{10}y_2
\end{equation}

\begin{equation}
    \frac{dy_{15}}{dt} = k_{14}y_1y_6
\end{equation}

\begin{equation}
    \frac{dy_{16}}{dt} = -k_{18}y_{16}-k_{19}y_{16}+k_{16}y_{4}
\end{equation}

\begin{equation}
    \frac{dy_{17}}{dt} = -k_{20}y_{17}y_6
\end{equation}

\begin{equation}
    \frac{dy_{18}}{dt} = k_{20}y_{17}y_6
\end{equation}

\begin{equation}
    \frac{dy_{19}}{dt} = -k_{21}y_{19}-k_{22}y_{19}-k_{24}y_{19}y_1+k_{23}y_1y_4+k_{25}y_{20}
\end{equation}

\begin{equation}
    \frac{dy_{20}}{dt} = -k_{25}y_{20}+k_{24}y_{19}y_1
\end{equation}

The initial conditions are 

\begin{gather*}
y = [0.0, 0.2, 0.0, 0.04, 0.0, 0.0, 0.1, 0.3, \\0.017, 0.0, 0.0, 0.0, 0.0, 0.0, 0.0, 0.0, 0.007, 0.0, 0.0, 0.0]\\ 
k = [0.25, 26.6, 12300.0, 0.00086, 0.00082, \\15000.0, 0.00013, 24000.0, 16500.0, 9000.0,\\0.022, 12000.0, 1.88, 16300.0, 4.8 \cdot 10^{6}, \\0.00035, 0.0175, 1.0 \cdot 10^{8}, 4.44 \cdot 10^{11}, 1240.0, \\2.1, 5.78, 0.0474, 1780.0, 3.12]\\
\end{gather*} 
The time span for integration is $t = (0.0 s, 60 s)$ \cite{34}. 

The comparison with several algorithms is done with relative tolerances in the range $(10^{-4}, 10^{-9})$ and corresponding absolute tolerances in the range $(10^{-8}, 10^{-13})$ against a reference solution with relative and absolute tolerances $10^{-14}$. 

\end{document}